# IMPLICIT EULER TIME DISCRETIZATION AND FDM WITH NEWTON METHOD IN NONLINEAR HEAT TRANSFER MODELING

Ph.D. Filipov S.[1], Prof. D.Sc. Faragó I.[2]

[1] Department of Computer Science, University of Chemical Technology and Metallurgy, Bulgaria
[2] Department of Applied Analysis and Computational Mathematics, MTA-ELTE Research Group, Eötvös Loránd University, Hungary
filipovstefan@yahoo.com, faragois@cs.elte.hu

**Abstract:** This paper considers one-dimensional heat transfer in a media with temperature-dependent thermal conductivity. To model the transient behavior of the system, we solve numerically the one-dimensional unsteady heat conduction equation with certain initial and boundary conditions. Contrary to the traditional approach, when the equation is first discretized in space and then in time, we first discretize the equation in time, whereby a sequence of nonlinear two-point boundary value problems is obtained. To carry out the time-discretization, we use the implicit Euler scheme. The second spatial derivative of the temperature is a nonlinear function of the temperature and the temperature gradient. We derive expressions for the partial derivatives of this nonlinear function. They are needed for the implementation of the Newton method. Then, we apply the finite difference method and solve the obtained nonlinear systems by Newton method. The approach is tested on real physical data for the dependence of the thermal conductivity on temperature in semiconductors. A MATLAB code is presented.

## 1. Introduction

We consider the one-dimensional unsteady heat conduction equation [1-3]

$$\rho c_p \frac{\partial u}{\partial t} = \frac{\partial}{\partial x}\left(\kappa(u)\frac{\partial u}{\partial x}\right), \tag{1}$$

where $u(x,t)$ is the temperature at position $x$ and time $t$, $\rho$ is the density, $c_p$ is the heat capacity at constant pressure, and $\kappa$ is the thermal conductivity of the media. We assume that $\rho$ and $c_p$ have constant values, but $\kappa$ depends on the temperature $u$. This assumption is often justifiable for certain temperature range (e.g. for silicon [4]). Performing the differentiation in the right-hand side of (1) we get

$$\rho c_p \frac{\partial u}{\partial t} = \partial_u \kappa(u)\left(\frac{\partial u}{\partial x}\right)^2 + \kappa(u)\frac{\partial^2 u}{\partial x^2}. \tag{2}$$

When $\kappa$ does not depend on $u$, i.e. $\partial_u \kappa(u) = 0$, then (2) is a linear (parabolic) partial differential equation. When $\partial_u \kappa(u) \neq 0$, then (2) is nonlinear. Equation (2) will be solved on the spatial interval $[a,b]$ subject to certain boundary and initial conditions:

$$u(a,t) = \alpha(t), u(b,t) = \beta(t), t > 0, \tag{3}$$
$$u(x,0) = u_0(x), x \in [a,b]. \tag{4}$$

The boundary conditions (3) give the temperature at the two ends as function of time. The initial condition (4) specifies the initial spatial distribution of the temperature.

## 2. Implicit Euler time discretization

For the linear problem ($\partial_u \kappa(u) = 0$), the numerical methods for solving (2)-(4) are well elaborated (finite difference, finite element methods, etc.). Usually, (2) is first discretized in space, whereby an initial-value (Cauchy) problem for first order ODE system is obtained. If the explicit Euler method is used to solve the Cauchy problem, then the method is stable only for $0 < D\tau/h^2 \leq 1/2$, where $D = \kappa/(\rho c_p)$ is the thermal diffusivity, $h$ is the discretization step in space, and $\tau$ is the discretization step in time.

Our approach is different. We first discretize (2) in time. Using a time-step $\tau$, the time line $t \geq 0$ is partition by equally separated mesh-points:

$$t_n = n\tau, \ n = 0, 1, 2, \ldots \tag{5}$$

Then, using implicit Euler scheme [5], equation (2) is discretized on the mesh (5):

$$\rho c_p \frac{u_n - u_{n-1}}{\tau} = \partial_u \kappa(u_n)\left(\frac{du_n}{dx}\right)^2 + \kappa(u_n)\frac{d^2 u_n}{dx^2}, \tag{6}$$

where $u_n = u_n(x)$ and $u_{n-1} = u_{n-1}(x)$ approximate the values of $u(x,t_n)$ and $u(x,t_{n-1})$, respectively. Equation (6) approximates the partial differential equation (2). The error is $O(\tau)$, hence the discretization scheme is first-order accurate in time. The method is stable, unlike the explicit method, i.e. evaluating the right-hand side of (2) at $u_{n-1}$, which is only conditionally stable.

Solving (6) for the second spatial derivative of the temperature $u_n$ we get

$$\frac{d^2 u_n}{dx^2} = f(u_n, v_n; u_{n-1}), \tag{7}$$

where $v_n = du_n/dx$ is the temperature gradient and $f$ is the following nonlinear function:

$$f(u_n, v_n; u_{n-1}) = \frac{\phi(u_n, v_n; u_{n-1})}{\kappa(u_n)}, \tag{8}$$

$$\phi(u_n, v_n; u_{n-1}) = \rho c_p \frac{u_n - u_{n-1}}{\tau} - \partial_u \kappa(u_n) v_n^2. \tag{9}$$

Equation (7), together with the boundary conditions $u_n(a) = \alpha(t_n)$, $u_n(b) = \beta(t_n)$, constitutes a nonlinear two-point boundary value problem (TPBVP) for the unknown function $u_n$. If the function $u_{n-1}$ is known (given) the problem can be solved by using some numerical technique for nonlinear problems. Thus, starting from the initial condition $u_0$ we can solve successively (7) for $n = 1, 2, \ldots$

## 3. Derivatives for the Newton method

The implementation of the Newton method requires the partial derivatives of $f(u_n, v_n; u_{n-1})$ with respect to $u_n$ and $v_n$. Introducing the notation $f_n = f(u_n, v_n; u_{n-1})$, $\phi_n = \phi(u_n, v_n; u_{n-1})$ and denoting the derivatives by $q_n = q(u_n, v_n; u_{n-1})$, $p_n = p(u_n, v_n)$ we get:

$$q_n = \frac{\partial f_n}{\partial u_n} = \frac{1}{\kappa(u_n)}\left(-f_n \partial_u \kappa(u_n) + \frac{\partial \phi_n}{\partial u_n}\right), \quad (10)$$

$$p_n = \frac{\partial f_n}{\partial v_n} = \frac{1}{\kappa(u_n)}\frac{\partial \phi_n}{\partial v_n}, \quad (11)$$

where

$$\frac{\partial \phi_n}{\partial u_n} = \frac{\rho c_p}{\tau} - \partial_{uu}^2 \kappa(u_n) v_n^2, \quad (12)$$

$$\frac{\partial \phi_n}{\partial v_n} = -2 \partial_u \kappa(u_n) v_n. \quad (13)$$

## 4. Finite difference method

Since $\rho c_p / \tau$ grows to infinity as $\tau$ goes to zero (which effects IVP solutions), it turns out that the finite difference method (FDM) [6] is a better choice for the solution of the obtained TPBVPs than the shooting method [6,7]. Hence, we adopt the FDM. The interval $[a, b]$ is partitioned by $N$ equally separated mesh-points:

$$x_i = a + (i-1)h, i = 1, 2, \ldots, N, h = \frac{b-a}{N-1}. \quad (14)$$

Equation (7) is discretizes on the uniform mesh (14) using the FDM with the central difference approximation:

$$\frac{u_{n,i+1} - 2u_{n,i} + u_{n,i-1}}{h^2} = f(u_{n,i}, v_{n,i}; u_{n-1,i}), \quad (15)$$

$$i = 2, 3, \ldots, N-1. \quad (16)$$

Correspondingly, everywhere in equations (8)-(13), we set $x = x_i$, and then replace the values of $u_n(x_i)$, $v_n(x_i)$, and $u_{n-1}(x_i)$ with their approximations $u_{n,i}$, $v_{n,i}$, and $u_{n-1,i}$ where

$$v_{n,i} = \frac{u_{n,i+1} - u_{n,i-1}}{2h}. \quad (17)$$

Equation (15) approximates (7) with error $O(h^2)$, i.e. it is second-order accurate in space. Equation (15) holds for the inner mesh-points. At the boundaries we apply the boundary conditions:

$$u_{n,1} = \alpha(t_n), u_{n,N} = \beta(t_n) \quad (18)$$

## 5. Solving the nonlinear system by Newton method

Introducing the column-vector $\mathbf{G}_n = [G_{n,1}, G_{n,2}, \ldots, G_{n,N}]^T$ with components

$$G_{n,1} = u_{n,1} - u_a(t_n), G_{n,N} = u_{n,N} - u_b(t_n), \quad (19)$$

$$G_{n,i} = u_{n,i+1} - 2u_{n,i} + u_{n,i-1} - h^2 f_{n,i}, \quad (20)$$

$$f_{n,i} = f(u_{n,i}, v_{n,i}; u_{n-1,i}), \quad (21)$$

the system of nonlinear equations (15) and the boundary conditions (18) are written as one equation:

$$\mathbf{G}_n(\mathbf{u}_n) = 0, \quad (22)$$

where

$$\mathbf{u}_n = [u_{n,1}, u_{n,2}, \ldots, u_{n,N}]^T. \quad (23)$$

Starting by some initial guess $\mathbf{u}_n^{(0)}$, the nonlinear system (22) can be solved by the Newton iterative method:

$$\mathbf{u}_n^{(k+1)} = \mathbf{u}_n^{(k)} - \left(\mathbf{L}_n^{(k)}\right)^{-1} \mathbf{G}_n\left(\mathbf{u}_n^{(k)}\right), \quad k = 0, 1, 2, \ldots \quad (24)$$

where $\mathbf{L}_n^{(k)}$ is the Jacobian of $\mathbf{G}_n$ with respect to $\mathbf{u}_n$ evaluated at $\mathbf{u}_n^{(k)}$:

$$\mathbf{L}_n^{(k)} = \frac{\partial \mathbf{G}_n}{\partial \mathbf{u}_n}\left(\mathbf{u}_n^{(k)}\right). \quad (25)$$

Calculating the elements of the Jacobian we get:

$$L_{n\,(1,1)}^{(k)} = 1, L_{n\,(N,N)}^{(k)} = 1, L_{n\,(i,i)}^{(k)} = -2 - h^2 q_{n,i}^{(k)}, \quad (26)$$

$$L_{n\,(i,i-1)}^{(k)} = 1 + \tfrac{1}{2} h p_{n,i}^{(k)}, L_{n\,(i,i+1)}^{(k)} = 1 - \tfrac{1}{2} h p_{n,i}^{(k)}, \quad (27)$$

$$q_{n,i}^{(k)} = q\left(u_{n,i}^{(k)}, v_{n,i}^{(k)}; u_{n-1,i}\right), p_{n,i}^{(k)} = p\left(u_{n,i}^{(k)}, v_{n,i}^{(k)}\right). \quad (28)$$

Iteration (24) is a one-step (two-level) iteration. Starting from some initial guess $\mathbf{u}_n^{(0)}$, we can find each next approximation $\mathbf{u}_n^{(k+1)}$, $k = 0, 1, 2, \ldots$ using (24). If the sequence is convergent, then the limiting vector $\mathbf{u}_n = \lim_{k \to \infty}\left(\mathbf{u}_n^{(k+1)}\right)$ is a solution to the nonlinear system (22). In practice, the iteration process is usually ended when

$$\left\| \mathbf{u}_n^{(k+1)} - \mathbf{u}_n^{(k)} \right\| < \epsilon. \quad (29)$$

This inequality is called a stopping criteria. The vector $\mathbf{u}_n^{(k+1)}$ is taken as approximate solution to (22). As an initial guess $\mathbf{u}_n^{(0)}$, we can use the solution $\mathbf{u}_{n-1}$ found at the previous step.

## 6. Computer experiment

Consider a thin homogenous rod, along the $x$-axis between the points $x = 1$ and $x = 3$, without heat sources and without radiation. The density $\rho$ and the heat capacity $c_p$ are constant, but the thermal conductivity $\kappa$ depends on the temperature as

$$\kappa = \kappa_0 \exp(\chi u). \quad (30)$$

Such a temperature dependence actually occurs in real physical systems, e.g. for silicon [4]. We choose the following values of the constants: $\rho = 1$, $c_p = 1$, $\kappa_0 = 0.1$. The temperature at the two ends is kept constant:

$$u(1, t) = 2, u(3, t) = 1, t > 0. \quad (31)$$

The initial temperature profile is

$$u(x, 0) = 2 - \frac{x-1}{2} + (x-1)(x-3), x \in [1,3]. \quad (32)$$

To find the time evolution of (32), we solve the partial differential equation (1) with boundary conditions (31) and initial conditions (32) by the method described in this paper. The equation is solved for $\chi = -1.0, -0.5, 0, 0.5, 1.0$. The step-size is chosen to be $\tau = 0.5$ with integration range $0 \leq t \leq 15$. The spatial interval $x \in [1,3]$ is discretized by $N = 41$ mesh-points, i.e. $h = 0.05$. The results are shown in Fig. 1

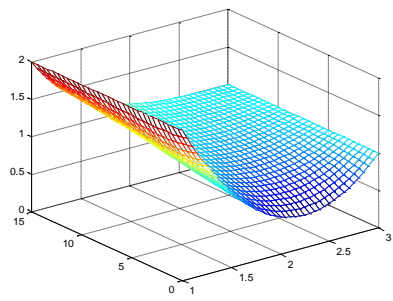
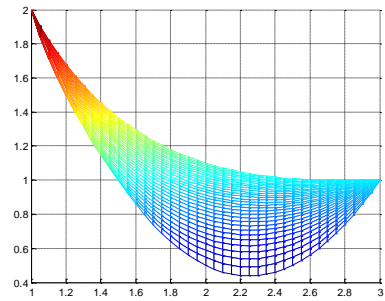

$\chi = -1.0$

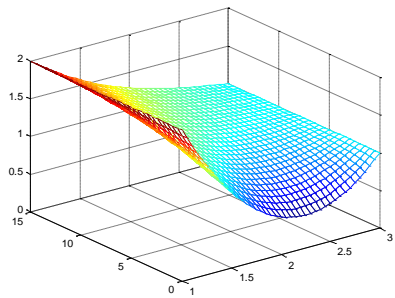
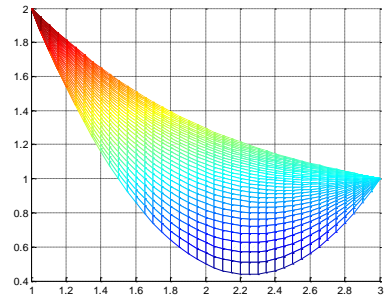

$\chi = -0.5$

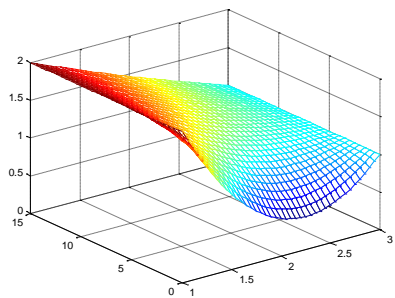
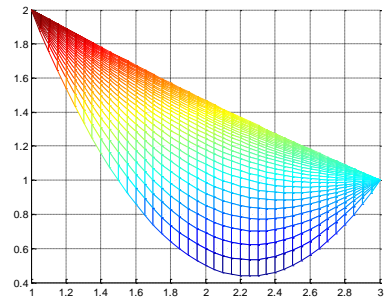

$\chi = 0$

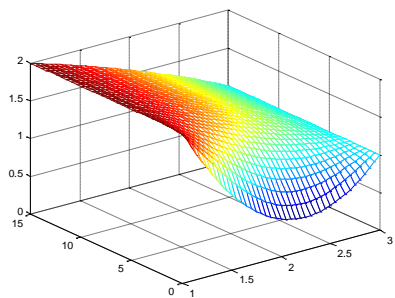
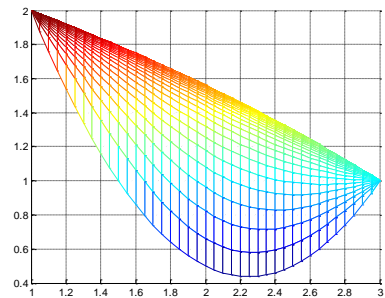

$\chi = 0.5$

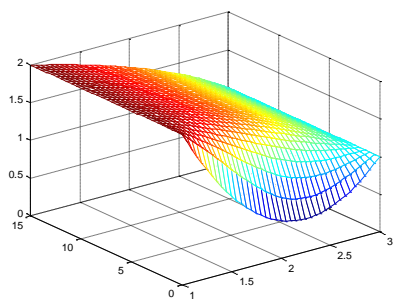
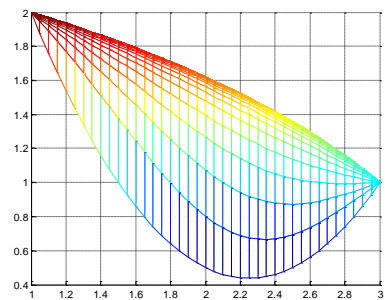

$\chi = 1.0$

***Fig. 1.*** *Solving the heat conduction equation (1) for boundary conditions (31), initial condition (32), and thermal conductivity (30). In addition to the 3D view (left), a side-view u vs. x is also shown (right).*

For χ = 0, 0.5, and 1.0 the final temperature distribution reached in the experiment is practically the steady-state distribution. For χ = −1.0 and − 0.5 a little bit more time is needed. The steady-state distribution for χ = 0 is, as expected, linear.

## 7. Conclusion

This paper considered heat transfer with temperature-dependent thermal conductivity. The one-dimensional unsteady heat conduction equation was solved numerically by using implicit time-discretization and FDM with Newton method for the solution of the arising nonlinear two-point boundary value problems. Data for the dependence of the thermal conductivity on temperature in certain semiconductors was used. The results obtained by the numerical computer experiments are consistent with the expected outcome. The proposed method is stable, unlike its explicit counterpart.

## 8. Appandix

A MATLAB code is presented for the numerical solution of the example provided in section 6.

```
function main
  rho=1; Cp=1; kappa0=0.1; chi=0.5;
  M=31; N=41;
  tEnd=15; tau=tEnd/(M-1);
  a=1; b=3; h=(b-a)/(N-1);
  alpha=2; beta=1;
  x=zeros(N,1);
  u0=zeros(N,1);
  for i=1:N
    x(i)=a+(i-1)*h;
    u0(i)=2-(x(i)-1)/2+(x(i)-1)*(x(i)-3);
  end
  t=zeros(M,1);
  for n=1:M
    t(n)=(n-1)*tau;
  end
  u_1=zeros(N,1); u=zeros(N,1);
  uNext=zeros(N,1); G=zeros(N,1);
  L=zeros(N,N); L(1,1)=1; L(N,N)=1;
  U=zeros(N,M); U(:,1)=u0;
  for n=2:M
    u=U(:,n-1); u_1=U(:,n-1);
    eps=1;
    while(eps>0.0001)
      G(1)=u(1)-alpha;
      G(N)=u(N)-beta;
      for i=2:N-1
        k=kappa0*exp(chi*u(i));
        Dk=chi*k;
        D2k=chi*Dk;
        v=(u(i+1)-u(i-1))/(2*h);
        A=rho*Cp/tau;
        phi=A*(u(i)-u_1(i))-Dk*v*v;
        f=phi/k;
        q=(-f*Dk+A-D2k*v*v)/k;
        p=-2*Dk*v/k;
        G(i)=u(i+1)-2*u(i)+u(i-1)-h*h*f;
        L(i,i-1)=1+0.5*h*p;
        L(i,i)=-2-h*h*q;
        L(i,i+1)=1-0.5*h*p;
      end
      uNext=u-L\G;
      eps=sqrt(h*(uNext-u)'*(uNext-u));
      u=uNext;
    end
    U(:,n)=u;
  end
  mesh(x,t,U');
end
```

The mathematical quantities and the corresponding variables used in the MATLAB code are shown in Table 1.

| Paper | MATLAB |
|---|---|
| $\rho$ | rho |
| $c_p$ | Cp |
| $\kappa_0$ | kappa0 |
| $\chi$ | chi |
| $\tau$ | tau |
| $N$ | N |
| $a$ | a |
| $b$ | b |
| $h$ | h |
| $\alpha$ | alpha |
| $\beta$ | beta |
| $x_i$ | x(i) |
| $t_n$ | t(n+1) |
| $\mathbf{u}_0$ | u0 |
| $\mathbf{u}_n^{(k)}$ | u |
| $\mathbf{u}_n^{(k+1)}$ | uNext |
| $\mathbf{u}_{n-1}$ | u_1 |
| $u_{n,i}$ | U(i,n+1) |
| $\mathbf{G}_n(\mathbf{u}_n^{(k)})$ | G |
| $\mathbf{L}_n^{(k)}$ | L |
| $\|\|\mathbf{u}_n^{(k+1)} - \mathbf{u}_n^{(k)}\|\|$ | eps |
| $u_{n,i}^{(k)}$ | u(i) |
| $\kappa(u_{n,i}^{(k)})$ | k |
| $\partial_u \kappa(u_{n,i}^{(k)})$ | Dk |
| $\partial_{uu}^2 \kappa(u_{n,i}^{(k)})$ | D2k |
| $v_{n,i}^{(k)}$ | v |
| $\phi(u_{n,i}^{(k)}, v_{n,i}^{(k)}; u_{n-1,i})$ | phi |
| $f(u_{n,i}^{(k)}, v_{n,i}^{(k)}; u_{n-1,i})$ | f |
| $q_{n,i}^{(k)}$ | q |
| $p_{n,i}^{(k)}$ | p |

*Table 1.*

## 9. Acknowledgement

This research was carried out in the ELTE Institutional Excellence Program (1783-3/2018/FEKUTSRAT) supported by the Hungarian Ministry of Human Capacities, and supported by the Hungarian Scientific Research Fund OTKA, No. K112157 and SNN125119

## 10. References

[1] R.B. Bird, W.E. Stewart, E.N. Lightfoot, Transport Phenomena, 2nd Edition, John Wiley & Sons, Inc., (2002)
[2] H.S. Carslaw, J.C. Jaeger, Conduction of Heat in Solids, 2nd Edition, Oxford University Press, (1986)
[3] S.L. Sobolev, Partial Differential Equations of Mathematical Physics, Dover Publications, (1989)
[4] J. Lienemann, A. Yousefi, J.G. Korvink, Nonlinear Heat Transfer Modeling, In: P. Benner, D.C. Sorensen, V. Mehrmann (eds) Dimension Reduction of Large-Scale Systems, Lecture Notes in Computational Science and Engineering, 45 (2005) 327-331
[5] U.M. Ascher, S.J. Ruuth, R.J. Spiteri: Implicit-Explicit Runge-Kutta Methods for Time-Dependent Partial Differential Equations, Appl Numer Math, vol. 25 (2-3) (1997)
[6] U.M. Ascher, R.M.M. Mattjeij, R.D. Russel, Numerical Solution of Boundary Value Problems for Ordinary Differential Equations, in: Classics in Applied Mathematics, vol. 13, SIAM, (1995)
[7] S.M. Filipov, I.D. Gospodinov, I. Faragó, Shooting-projection method for two-point boundary value problems, Appl. Math. Lett. 72 (2017) 10-15
https://doi.org/10.1016/j.aml.2017.04.002